\documentclass{amsproc}

\newtheorem{theorem}{Theorem}[section]
\newtheorem{lemma}[theorem]{Lemma}

\theoremstyle{definition}

\newtheorem{proposition}[theorem]{Proposition}
\newtheorem{corollary}[theorem]{Corollary}

\theoremstyle{remark}

\newtheorem{definition}[theorem]{Definition}
\newtheorem*{remark*}{Remark}

\newtheorem*{remarks}{Remarks}

\numberwithin{equation}{section}

\usepackage{amssymb}
\usepackage{amsfonts}
\usepackage{latexsym}

\newcommand{\Ext}{\text{Ext}}
\newcommand{\Hom}{\text{Hom}}

\newcommand{\ben}{\begin{enumerate}}
\newcommand{\een}{\end{enumerate}}

\newcommand{\CC}{{\mathbb{C}}}
\newcommand{\ZZ}{{\mathbb{Z}}}

\begin{document}

\title{Quantization, Orbifold Cohomology, and Cherednik Algebras}

\author{Pavel Etingof}
\address{Department of Mathematics, Massachusetts Institute of Technology,
Cambridge, MA 02139, USA}
\email{etingof@math.mit.edu}
\thanks{The work of the first author
was supported by the NSF grant DMS-9988796.}

\author{Alexei Oblomkov}
\address{Department of Mathematics, Massachusetts Institute of Technology,
Cambridge, MA 02139, USA}
\email{oblomkov@math.mit.edu}


\begin{abstract}
We compute the Hochschild homology of the crossed product
$\Bbb C[S_n]\ltimes A^{\otimes n}$ in terms of the Hochschild homology
of the associative algebra $A$ (over $\Bbb C$). It allows us to compute
the Hochschild (co)homology of $\Bbb C[W]\ltimes A^{\otimes n}$ where $A$
is the $q$-Weyl algebra or any its degeneration and $W$ is the Weyl group
of type $A_{n-1}$ or $B_n$. For a deformation quantization $A_+$ of an
affine symplectic variety $X$ we show that the Hochschild homology of
$S^n A$, $A=A_+[\hbar^{-1}]$ is additively isomorphic to the Chen-Ruan
orbifold cohomology of $S^nX$ with coefficients in $\Bbb C((\hbar))$.
We prove that for $X$ satisfying $H^1(X,\Bbb C)=0$ (or $A\in VB(d)$)
the deformation of $S^nX$ ($\Bbb C[S_n]\ltimes A^{\otimes n}$)
which does not come from deformations of $X$ ($A$) exists
if and only if $\dim X=2$ ($d=2$). In particular if $A$ is
$q$-Weyl algebra (its trigonometric or rational degeneration)
then the corresponding nontrivial deformations yield the double
affine Hecke algebras of type $A_{n-1}$
(its trigonometric or rational versions) introduced by Cherednik.
\end{abstract}

\maketitle

\section{Introduction}
In this note, given an associative algebra $A$ over $\Bbb C$, we compute
the Hochschild homology of the crossed product
$\Bbb C[S_n]\ltimes A^{\otimes n}$. If $A$ is simple,
$\dim(A)=\infty$, and the center $Z(A)$ of $A$ is $\Bbb C$, then this homology
coincides with the homology of $S^nA$.

If $A$ satisfies the ``Gorenstein'' properties of \cite{VB1,VB2},
then this computation allows us to compute the Hochschild
cohomology of $\Bbb C[S_n]\ltimes A^{\otimes n}$, and of $S^nA$ for
simple $A$ with $\dim(A)=\infty$ and $Z(A)=\Bbb C$.

In particular, we obtain a result conjectured
(in a much stronger form) by Ginzburg and Kaledin (\cite{GK}, (1.3)):
if $X$ is an affine symplectic
algebraic variety over $\Bbb C$, $A_+$ is a deformation quantization
of $X$, and $A=A_+[\hbar^{-1}]$,
then the Hochschild cohomology of the algebra $S^nA$
(which is a quantization of the singular Poisson variety $S^nX=X^n/S_n$)
is additively isomorphic to the Chen-Ruan orbifold (=stringy) cohomology
of $S^nX$ with coefficients in $\Bbb C((\hbar))$.
\footnote{It is useful to compare
this result with the following theorem of Brylinski-Kontsevich:
if $Y$ is an affine symplectic variety,
$B_+$ a deformation quantization of $Y$, and
$B=B_+[\hbar^{-1}]$,
then the Hochschild cohomology
$HH^*(B)$ is isomorphic to the topological (=singular) cohomology
$H^*(Y,\Bbb C((\hbar)))$. The main difference between the two
cases is that $Y$ is smooth, while $S^nX$ is not.}
If $X$ is a surface, this cohomology is isomorphic
to the cohomology of the Hilbert scheme
${\rm Hilb}_n(X)$ (the G\"ottsche formula, \cite{Go}).

As a corollary, we get that if $H^1(X,\Bbb C)=0$ then
for $n>1$, $\dim HH^2(S^nA)=\dim HH^2(A)$ if $\dim X>2$, and
$\dim HH^2(S^nA)=\dim HH^2(A)+1$ if $\dim X=2$ (i.e., $X$ is a
surface). This implies that if $\dim X>2$ then all deformations
of $S^nA$ come from deformations of $A$; on the other hand,
if $X$ is a surface, then there is an
additional parameter of deformations of $S^nA$, which does not
come from deformations of $A$. In this case, the universal deformation
of $S^nA$ is a very interesting algebra.
For example, if $X=\Bbb C^2$, it is the spherical subalgebra of
the rational Cherednik algebra
attached to the group $S_n$ and its standard representation
$\Bbb C^n$ (see e.g. \cite{EG}).
In general, this deformation exists only over formal
series, but we expect that in ``good'' cases (when $X$ has a
compactification to which the Poisson bracket extends by zero at infinity,
see \cite{Ko2},\cite{Ar}), this deformation exists ``nonperturbatively''.

We note that many of the results below are apparently known to
experts; we present them with proofs since we could not find an
exposition of them in the literature.

{\bf Acknowledgments}
The authors thank M. Artin, Yu. Berest, P. Bressler,
V. Dolgushev, V. Ginzburg, L. Hesselholt,
D. Kaledin, and M. Kontsevich for useful discussions.

\section{Homology with twisted coefficients}

Let $A$ be an associative unital algebra over $\Bbb C$.
Let $M$ be an $A$-bimodule. Let $\sigma$ be the cyclic
permutation (12...n). Denote by $(A^{\otimes (n-1)}\otimes M)\sigma$
the space $A^{\otimes (n-1)}\otimes M$ with the
structure of $A^{\otimes n}$-bimodule,
given by the formula
$$
(a_1\otimes...\otimes a_n)(b_1\otimes...\otimes b_{n-1}\otimes m)(c_1\otimes...\otimes c_n)=
a_1b_1c_2\otimes...\otimes a_{n-1}b_{n-1}c_n\otimes a_nmc_1.
$$

\begin{proposition}\label{homolog}
(i) There exists a natural isomorphism of Hochschild homology
$HH_i(A^{\otimes n},(A^{\otimes (n-1)}\otimes M)\sigma)\to
HH_i(A,M)$.

(ii) $\sigma$ acts trivially on
$HH_i(A^{\otimes n},A^{\otimes n}\sigma)$.
\end{proposition}

\begin{proof}
(i) Recall that $HH_i(A,M)=Tor_i(A,M)$, the derived functors
of $Tor_0(A,M)=A\otimes_{A\otimes \bar A}M=M/[A,M]$ (here $\bar A$
is the opposite algebra of $A$). These spaces are computed as follows.
Let $...F_1\to F_0\to M$ be a free resolution of $M$ as an
$A\otimes \bar A$-module. Then we have a complex
$...F_1/[A,F_1]\to F_0/[A,F_0]\to 0$, and it's
i-th homology is $HH_i(A,M)$.

Now let us compute $HH_i(A^{\otimes n},(A^{\otimes (n-1)}\otimes M)\sigma)$.
The above resolution yields a resolution
$$
...(A^{\otimes (n-1)}\otimes F_1)\sigma\to
(A^{\otimes (n-1)}\otimes F_0)\sigma
\to  (A^{\otimes (n-1)}\otimes M)\sigma
$$
This resolution is not free (since its terms
are not free $A^{\otimes n}\otimes \bar A^{\otimes n}$-modules),
but nevertheless it can be used to compute
the homology groups \linebreak
$Tor_i(A^{\otimes n},(A^{\otimes (n-1)}\otimes M)\sigma)$.
This is a consequence of the following lemma.

\begin{lemma}
Let $F$ is a free $A\otimes \bar A$-module. Then for $i>0$ we have
$$
Tor_i(A^{\otimes n},(A^{\otimes (n-1)}\otimes F)\sigma)=0.
$$
\end{lemma}

\begin{proof}
It is sufficient to consider the case when
$F=A\otimes A$, with the bimodule structure $a(b\otimes c)d=ab\otimes cd$.
In this case
$(A^{\otimes (n-1)}\otimes F)\sigma=A^{\otimes (n+1)}$,
with the $A^{\otimes n}$-bimodule structure given by the formula
$$
(a_1\otimes...\otimes a_n)(b_1\otimes...\otimes b_{n+1})(c_1\otimes...\otimes c_n)=
a_1b_1c_2\otimes...\otimes a_{n-1}b_{n-1}c_n\otimes a_nb_n\otimes b_{n+1}c_1
$$
This shows that the $A^{\otimes n}\otimes \bar A^{\otimes n}$
module $A^{\otimes (n+1)}$ is induced (by adding the \linebreak
$n+1$-th component) from the module $Y=A^{\otimes n}$ over the
subalgebra $B=A^{\otimes n}\otimes \bar A^{\otimes (n-1)}$ spanned
by elements $a_1\otimes...\otimes a_n\otimes 1\otimes
c_2\otimes...\otimes c_n$. Therefore, by the Shapiro lemma,
$Tor_i(A^{\otimes n},(A^{\otimes (n-1)}\otimes F)\sigma)=
Tor_i(A^{\otimes n},Y)$ (where on the right hand side the modules
are over the algebra $B$). Now, the module $A^{\otimes n}$, as a
B-module, is also induced. Namely, it is induced (by adding the
first component) from the module $A^{\otimes (n-1)}$ over the
subalgebra $A^{\otimes (n-1)}\otimes \bar A^{\otimes (n-1)}$
spanned by elements \linebreak $1\otimes a_2\otimes...\otimes
a_n\otimes 1\otimes c_2\otimes...\otimes c_n$. Thus we can use the
Shapiro lemma again, and get $Tor_i(A^{\otimes n},(A^{\otimes
(n-1)}\otimes F)\sigma)= Tor_i(A^{\otimes (n-1)},A^{\otimes n})$,
where the $A^{\otimes (n-1)}$-bimodule structure on $A^{\otimes
n}$ is given by the formula
$$
(a_2\otimes...\otimes a_n)(b_1\otimes...\otimes b_{n})(c_2\otimes...\otimes c_n)=
b_1c_2\otimes a_2b_2c_3\otimes...\otimes a_{n-1}b_{n-1}c_n\otimes a_nb_n.
$$
Continuing to apply the
Shapiro lemma as above, we will eventually reduce the algebra
over which the Tor functors are computed to the ground field $\Bbb C$.
This implies that the higher Tor functors vanish, and the lemma is proved.
\end{proof}

Thus, to prove part (i) the proposition, it suffices to show
$$
A^{\otimes n}\otimes_{A^{\otimes n}\otimes \bar A^{\otimes n}}
(A^{\otimes (n-1)}\otimes M)\sigma=M/[A,M].
$$
This is straightforward.

(ii) Let $C$ be an $n\times m$-matrix. Let $s(C)$ denote the new
matrix, obtained from $C$ by cyclically permuting the columns
(putting the first column at the end), and then applying $\sigma$
to this column.

If $C$ is a matrix of elements of $A$, then $C$ can be regarded
as an $m-1$ chain of $A^{\otimes n}$ with coefficients
in $A^{\otimes n}\sigma$, by taking the tensor product
$$
c_{11}\otimes c_{21}\otimes...\otimes c_{n1}\otimes
c_{12}...\otimes c_{n2}\otimes ...\otimes c_{1m}\otimes...\otimes
c_{nm}.
$$
Thus, $s$ defines a linear operator on $m-1$-chains.

Let $C$ be an $m-1$-cycle of $A^{\otimes n}$ with coefficients
in $A^{\otimes n}\sigma$. We claim that
$$
C-\sigma(C)=
d\left(\sum_{j=0}^{m-1}(-1)^{j(m-1)}s^j(C)\otimes 1^{\otimes n}\right).
$$
This is checked by an easy direct computation
and implies the statement.
\end{proof}

\begin{remark*} L. Hesselholt explained to us that Proposition
\ref{homolog} is known in algebraic topology. A topological proof
of this proposition (see \cite{BHM}, Section 1) is as follows.
The Hochschild complex of $A$ with coefficients in $M$
can be viewed as a simplicial abelian group $X$; its homotopy groups
are the Hochschild homology groups.
The Hochschild complex of $A^{\otimes n}$ with
coefficients in $(A^{\otimes n-1}\otimes M)\sigma$
is canonically isomorphic to the simplicial abelian group $sd_n(X)$ (the edgewise
subdivision of $X$). Then part (i) of Proposition \ref{homolog}
follows from the fact that $sd_n(X)$ is homeomorphic to $X$ for
any simplicial set $X$. Part (ii) follows from the fact that for
$M=A$ the simplicial abelian group $X$ is cyclic (in the sense of Connes),
and hence the natural action of $\Bbb Z_n$ on $sd_n(X)$
extends to a continuous circle action.
\end{remark*}

\section{Homology and cohomology of crossed products}

\subsection{Homology of crossed products}
Let $A(n):=\Bbb C[S_n]\ltimes A^{\otimes n}$.
Let $P_n$ be the set of partitions of $n$.
If $\lambda$ is a partition of $n$, we set $p_i(\lambda), i\ge 1$
to be the multiplicity of occurrence of $i$ in $\lambda$.
For any algebra $B$, for brevity we write $HH_i(B)$ and $HH^i(B)$ for
$HH_i(B,B)$ and $HH^i(B,B)$.

\begin{theorem}\label{homo}
We have
$$
HH_*(A(n))=\oplus_{\lambda\in P_n}\bigotimes_{i\ge 1}
S^{p_i(\lambda)}HH_*(A).
$$
Here $HH_*$ is regarded as a functor with values in the category
of supervector spaces, with $HH_i$ being of parity $(-1)^i$.
\end{theorem}

\begin{remark*} If $A$ is the Weyl algebra ${\mathcal A}_1$, this result
follows from \cite{AFLS}.
\end{remark*}

\begin{proof}
Let $G$ be a finite group acting on an algebra
$B$, and let $D=\Bbb C[G]\ltimes B$.

\begin{proposition}\label{afls} (see e.g. \cite{AFLS})
We have
$$
HH_i(D)=\oplus_{C\subset G}(\oplus_{g\in C}HH_i(B,Bg))^G
$$
where $C$ runs over conjugacy classes of $G$.
\end{proposition}

Let us apply this formula to the case $B=A^{\otimes n}$, $G=S_n$.
Recall that \linebreak $HH_*(B_1\otimes B_2)=
HH_*(B_1)\otimes HH_*(B_2)$, and
that conjugacy classes in $S_n$ are labeled by
partitions. These facts and Proposition \ref{homolog}
imply the result.
\end{proof}

\begin{corollary}\label{simp}
Suppose that $A$ is an infinite dimensional simple algebra with trivial center.
Then
$$
HH_*(S^nA)=\oplus_{\lambda\in P_n}\bigotimes_{i\ge 1}
S^{p_i(\lambda)}HH_*(A).
$$
\end{corollary}

\begin{proof} Since the center of $A$ is trivial, by \cite{MR}, Lemma 9.6.9,
the algebra $A^{\otimes n}$ is simple.

We claim that any element $\sigma\in S_n$, $\sigma\ne 1$, defines
an outer automorphism of $A^{\otimes n}$. Indeed,
by conjugating $\sigma$ we may assume that $\sigma(1)=m\ne 1$,
$\sigma^{-1}(1)=p\ne 1$. Assume the contrary, i.e. that
$\sigma(z)=gzg^{-1}$, $g\in A^{\otimes n}$.  Then
for any $a\in A$, $ga_1=a_mg, a_1g=ga_p$, where $a_i$ denotes
the tensor product of $a$ put in the $i$-th component with units
in the other components. This implies that the finite dimensional
space $I$ spanned by the first components of $g$ is a 2-sided
ideal in $A$, which is a contradiction with the facts that
$A$ is simple and infinite dimensional.

By \cite{NZ}
(see also \cite{Mo}), if $B$ is a simple algebra and $G$ a finite group
acting on $B$ by outer automorphisms then
the algebra $G\ltimes B$ is
simple. It follows that $A(n)$ is a simple algebra.
Hence the algebra $A(n)$ is Morita equivalent to $S^nA$.
This fact and Theorem \ref{homo}
implies the statement.
\end{proof}

\subsection{The class $VB(d)$ of algebras}
Let us now define a class of ``Gorenstein'' algebras $VB(d)$,
introduced by Van den Bergh in \cite{VB1,VB2}.

\begin{definition}
We will say that an algebra $A$ is in the class $VB(d)$
if $A$ has finite Hochschild dimension, and
$\Ext_{A\otimes \bar{A}}^*(A,A\otimes \bar A)$ is concentrated in
degree $d$, where it equals to $A$ (as an $A$-bimodule).
\end{definition}

For example, if $A={\mathcal O}_X$, where $X$ is a smooth affine
algebraic variety of dimension $d$, then $A\in VB(d)$
if and only if the canonical bundle of $X$ is trivial, i.e.
$X$ admits a nonvanishing differential form of top degree
(a volume form).

Let $G$ be a finite group acting by automorphisms
of an algebra $A\in VB(d)$. We say that the action of $G$ is
unimodular if there exists an isomorphism of bimodules
$\Ext_{A\otimes \bar{A}}^d(A,A\otimes \bar A)\to A$
which commutes with $G$. For example, in the case $A={\mathcal O}_X$,
where $X$ a smooth affine variety of dimension $d$, the action
of $G$ is unimodular if and only if there exists a $G$-invariant
volume form on $X$.

In the sequel we'll need the following proposition about the
class $VB(d)$.

\begin{proposition} \label{GA}
If $A\in VB(d)$ and $G$ is a finite group acting on
$A$ in a unimodular fashion then $\Bbb C[G]\ltimes A\in VB(d)$.
\end{proposition}

\begin{proof}
The proof is straightforward from Shapiro's lemma.
\end{proof}

\subsection{Cohomology of crossed products}

\begin{corollary}\label{vdb}
Let $A\in VB(d)$ for even $d$. Then
$$
HH^*(A(n))=\oplus_{\lambda\in P_n}\bigotimes_{i\ge 1}
S^{p_i(\lambda)}HH^{*-d(i-1)}(A),
$$
with $HH^i$ of parity $(-1)^i$.
\end{corollary}

\begin{proof} It is easy to see that if $d$ is even and $A\in VB(d)$ then
$A^{\otimes n}\in VB(nd)$, and the action of $S_n$ on $A^{\otimes n}$
is unimodular. Thus by Proposition \ref{GA}, $A(n)\in VB(nd)$. So
according to \cite{VB1,VB2}, $HH^*(A)=HH_{d-*}(A)$,
and $HH^*(A(n))=HH_{dn-*}(A(n))$.
This fact and Theorem \ref{homo} implies the statement.
\end{proof}

\begin{corollary} \label{vdb1}
If $d$ is even and $A\in VB(d)$ is simple infinite dimensional
with trivial center, then
$$
HH^*(S^nA)=\oplus_{\lambda\in P_n}\bigotimes_{i\ge 1}
S^{p_i(\lambda)}HH^{*-d(i-1)}(A).
$$
\end{corollary}

\begin{proof}
The same as for Corollary \ref{simp}.
\end{proof}

\begin{remark*} Let
$HH^*(A(\bullet)):=\oplus_{n\ge 0}HH^*(A(n))$.
Corollary \ref{vdb} implies that
$$
HH^*(A(\bullet))=\otimes_{j\ge 0}S^\bullet HH^{*-dj}(A)
$$
This implies the following formula
for the generating function
of $\dim HH^i(A(n))$: if $b_k$ are the Betti numbers of $A$,
then
\begin{equation}\label{Goe}
\sum t^iq^n\dim HH^i(A(n))=\prod_{m\ge 1}\prod_{k\ge 0}
(1+(-1)^{k-1}q^mt^{k+d(m-1)})^{(-1)^{k-1}b_k}.
\end{equation}
\end{remark*}

\section{Hochschild cohomology of quantizations}

Now we consider the following situation.
Let $X$ be an affine symplectic
algebraic variety over $\Bbb C$, $A_+$ be a deformation quantization
of $X$, and $A=A_+[\hbar^{-1}]$.
In this case the algebra $S^nA_+$ is a deformation quantization
of the singular Poisson variety $S^nX$, and $S^nA$
is its quantization over $\Bbb C((\hbar))$.

\begin{remark*} The algebra $A$ is over
$\Bbb C((\hbar))$ (rather than $\Bbb C$)
and carries an $\hbar$-adic
topology; this will not be important
to us, except that all tensor products will
have to be completed with respect to this topology.
\end{remark*}

\begin{theorem}\label{Hoh}
The Hochschild cohomology
$HH^*(S^nA)$ is (additively) isomorphic
to the orbifold cohomology
$H^*_{orb}(S^nX)\otimes\Bbb C((\hbar))$ \cite{Ba,CR}.
\end{theorem}

\begin{remark*} This is a very special case of the conjecture (1.3) from
\cite{GK}, which states that the above isomorphism is multiplicative and
extends from symmetric powers to any symplectic orbifold.
\end{remark*}

\begin{proof}
 Let us first show that $A\in VB(d)$, where $d=\dim X$.

1. By a deformation argument, the algebra $A$ is of finite Hochschild dimension
(since so is the structure ring $A_0=O(X)$).

Further, the cohomology
$\Ext_{A\otimes \bar{A}}^*(A,A\otimes \bar A)$ is concentrated in
degree $d$, since so is the corresponding cohomology for
$A_0$ instead of $A$.

Moreover, it is clear by deformation argument that
$U:=\Ext_{A\otimes \bar{A}}^d(A,A\otimes \bar A)$
is an invertible $A$-bimodule.
Thus, by \cite{VB1,VB2},
$HH^i(A,N)=HH_{d-i}(A,U\otimes_A N)$.

2. Let us show that $U=A$.

$U^{-1}$ is clearly free
as a left and right module over $A$ (since it is such for $A_0$).
Thus, $U^{-1}=A\gamma$, where $\gamma$ is an automorphism of $A_+$
which equals $1$ modulo $\hbar$.

We have $HH^0(A,U^{-1})=HH_d(A)$.
But $HH_d(A)$ equals the Poisson homology
$HP_d(A_0)\otimes \Bbb C((\hbar))$
(see e.g. \cite{NeTs}, Theorem A2.1;
the theorem is in the $C^\infty$-setting but the proof applies to
smooth affine algebraic varieties as well).
By Brylinski's theorem \cite{Br},
$HP_d(A_0)=H^0(X,\Bbb C)=\Bbb
C$. Thus,
$HH^0(A,U^{-1})\ne 0$,
and there exists an element $x\in A_+$ such that
$ax=x\gamma(a)$ for all $a\in A_+$.
Clearly, we can assume that $x$ projects to $x_0\ne 0$
in $A_0=A_+/\hbar A_+$.
Let $\gamma_1: A_0\to A_0$ be the reduction of
$(\gamma-{\rm id})/\hbar$ modulo $\hbar$. We find
$\lbrace{a,x_0\rbrace}=x_0\gamma_1(a)$, $a\in A_0$.
This implies that $x_0$ is nonvanishing on $X$.
So $x$ is invertible and $\gamma$ is inner.
Thus $U=A$, and $A\in VB(d)$.

3. Moreover, the algebra $A$ is simple.
Indeed, if $I$ is a nonzero two-sided ideal in $A$, then
$I_+:=I\cap A_+$ is a saturated ideal in $A_+$.
Let $I_0=I_+/\hbar I_+\subset O(X)=A_+/\hbar A_+$.
Then $I_0$ is a nonzero Poisson ideal in $O(X)$.
But $X$ is symplectic, so $I_0=O(X)$. Since $I_+$ is saturated,
$I_+=A_+$ and hence $I=A$.

Similarly, the algebra
$A^{\otimes n}$ is simple.
Also, the group $S_n$ acts on $A^{\otimes n}$ by outer
automorphisms (as $A^{\otimes n}$ has no nontrivial
inner automorphisms of finite order).
Thus, Corollaries \ref{simp} and \ref{vdb1} apply to
$A$.

4. Now, by Grothendieck's algebraic De Rham theorem,
the topological cohomology of $X$ is isomorphic
to its algebraic de Rham cohomology.
By Brylinski's theorem \cite{Br}, the latter is isomorphic
to the Poisson cohomology $HP^*(X)$.
By Kontsevich's formality theorem \cite{Ko1}
$HP^*(X)\otimes \Bbb C((\hbar))$ is isomorphic to $HH^*(A)$.
This and Corollary \ref{vdb1} imply that
\begin{equation}\label{Goet}
HH^*(S^nA)=\oplus_{\lambda\in P_n}\bigotimes_{i\ge 1}
S^{p_i(\lambda)}H^{*-d(i-1)}(X,\Bbb C((\hbar))).
\end{equation}
By \cite{CR},\cite{Uribe}, this is isomorphic to
$H^*_{orb}(S^nX)\otimes \Bbb C((\hbar))$, as desired.
\end{proof}

\begin{remark*} It is easy to see that in the case when $X$ is a
surface, formula (\ref{Goet}) turns into the G\"ottsche formula
\cite{Go} for the generating function of the Poincar\'e
polynomials of the Hilbert schemes $Hilb_n(X)$.
\end{remark*}

\section{Examples}
Now we will apply the above results to some particular algebras.
\subsection{Weyl algebras}

Let ${\mathcal A}$ be the Weyl algebra generated by $x,p$
with definiting relation
\begin{gather}
[x,p]=1.
\end{gather}

Let ${\Bbb A}$ be the trigonometric Weyl algebra
generated by $X^{\pm 1},p$ with defining relation
\begin{gather}
[X,p]=X.
\end{gather}

Let ${\mathbf A}$ be the q-Weyl algebra
generated by $X^{\pm 1}, P^{\pm 1}$
with defining relation
\begin{gather}
XP=qPX.
\end{gather}

We always suppose that $q$ is not a root of unity.

Let $\epsilon$ be the automorphism of ${\mathcal A}$
changing the sign of $x,p$. We also denote by $\epsilon$ the
automorphism of ${\Bbb A}$
sending $X$ to $X^{-1}$ and $p$ to $-p$, and
the automorphism of ${\mathbf A}$
that sends $X$ and $P$ to their inverses.

Using $\epsilon$,
one can define, in an obvious way,
 a natural action of the classical Weyl groups $W$ of type
$A_{n-1},B_n$ on the algebras $\mathcal{A}^{\otimes n}$
${\Bbb A}^{\otimes n}$, and
$\mathbf{A}^{\otimes n}$. Here we will give formulas (a la G\"ottsche
formula) for the generating series of the dimensions of
the Hochschild cohomology
of $\Bbb C[W]\ltimes \mathcal{A}^{\otimes n}$,
$\Bbb C[W]\ltimes {\Bbb A}^{\otimes n}$ and
$\Bbb C[W]\ltimes \mathbf{A}^{\otimes n}$. In the case of
$\Bbb C[S_n]\ltimes \mathcal{A}^{\otimes n}$ this formula was derived in
\cite{AFLS}.

\begin{remark*} Formulas for $W$ of type $D_n$ can be obtained in a
similar way, but are more complicated.
\end{remark*}

Define the generating functions
\begin{gather*}
\mathcal{P}_W=\sum t^iq^n \dim HH^i(\Bbb C[W]\ltimes\mathcal{A}^{\otimes n}),\\
{\Bbb P}_W=\sum t^iq^n \dim HH^i(\Bbb C[W]\ltimes {\Bbb A}^{\otimes n}),\\
\mathbf{P}_W=\sum t^iq^n \dim HH^i(\Bbb C[W]\ltimes\mathbf{A}^{\otimes n}),
\end{gather*}

\begin{theorem}\label{classicalgps} The generating functions
for $W$ of types $A,B$ are given by the following formulas:
\begin{gather*}
\mathcal{P}_{{A}}=\prod_{m\ge 1}(1-q^m t^{2(m-1)})^{-1},\\
{\Bbb P}_{{A}}=\prod_{m\ge 1}(1-q^m t^{2(m-1)})^{-1}(1+q^m t^{2m-1}),\\
\mathbf{P}_{{A}}=\prod_{m\ge 1}(1-q^mt^{2(m-1)})^{-1}(1-q^m
t^{2m})^{-1}(1+q^m t^{2m-1})^2,\\
\mathcal{P}_{{B}}=\prod_{m\ge 1}(1-q^m t^{2(m-1)})^{-1} (1-q^m
t^{2m})^{-1},\\
{\Bbb P}_{{B}}=\prod_{m\ge 1}(1-q^m t^{2(m-1)})^{-1} (1-q^m
t^{2m})^{-2},\\
\mathbf{P}_{{B}}=
\prod_{m\ge 1}(1-q^m t^{2(m-1)})^{-1}(1-q^m t^{2m})^{-5}.
\end{gather*}
\end{theorem}

The proof of this theorem is given below.

\begin{remark*} Let $\Gamma$ be a finite subgroup of $SL_2(\Bbb C)$.
Then the wreath product $W=\Bbb S_n\ltimes \Gamma^n$ acts
on ${\mathcal A}^{\otimes n}$. Similarly to the theorem, one can show that
the generating function ${\mathcal P}_\Gamma$ of dimensions
of the Hochschild cohomology of $\Bbb C [W]\ltimes A^{\otimes n}$ is
$$
\mathcal{P}_{\Gamma}=\prod_{m\ge 1}(1-q^m t^{2(m-1)})^{-1} (1-q^m
t^{2m})^{1-\nu(\Gamma)},
$$
where $\nu(\Gamma)$ is the number of conjugacy classes of $\Gamma$.
\end{remark*}

\subsection{Hochschild cohomology of ${\mathcal A}$,
${\Bbb A}$, and $\mathbf{A}$}
\label{homolAA}

\begin{proposition} \label{p1}
\begin{gather*}
HH^{i}(\mathcal{A})=0, \mbox{ if } i\ne 0,\\
HH^0(\mathcal{A})=\CC.
\end{gather*}
\end{proposition}

\begin{proposition} \label{p2}
\begin{gather*}
HH^{i}({\Bbb A})=0, \mbox{ if } i\ne 0,1\\
HH^0({\Bbb A})=\CC,\
HH^1({\Bbb A})=\CC.
\end{gather*}
\end{proposition}

\begin{proposition} \label{p3}
Suppose that $q$ is not a root of unity. Then
\begin{gather*}
HH^{i}(\mathbf{A})=0,\mbox{ if } i\ne 0,1,2,\\
HH^1(\mathbf{A})=\CC^2,\\
HH^0(\mathbf{A})=HH^2(\mathbf{A})=\CC.
\end{gather*}
\end{proposition}

These propositions are well known.
Proposition \ref{p1} is proved in \cite{Srid}.
For the proof of Proposition \ref{p3}
see e.g. \cite{Ob}, Section 5. Proposition \ref{p2}
is proved by a similar argument to \cite{Ob}.

\begin{proposition}\label{AxZ2}
(\cite{Ob}, Section 5). Suppose that $q$ is not a root of unity. Then
\begin{gather*}
HH^i(\CC [\ZZ_2]\ltimes \mathcal{A})=
HH^i(\CC [\ZZ_2]\ltimes {\Bbb A})=
HH^i(\CC [\ZZ_2]\ltimes \mathbf{A})=0,
\mbox{ if } i\ne
0,2,\\
HH^0(\CC [\ZZ_2]\ltimes \mathcal{A})=
HH^0(\CC [\ZZ_2]\ltimes {\Bbb A})=
HH^0(\CC [\ZZ_2]\ltimes \mathbf{A})=\Bbb C,
\\
HH^2(\CC [\ZZ_2]\ltimes \mathcal{A})=\Bbb C,
\\
HH^2(\CC [\ZZ_2]\ltimes {\Bbb A})=\Bbb C^2,
\\
HH^2(\CC [\ZZ_2]\ltimes \mathbf{A})=\Bbb C^5.
\end{gather*}
\end{proposition}

\subsection{Duality}

\begin{proposition} \label{VB(2)}
The algebras ${\mathcal A}$, ${\Bbb A}$, and ${\mathbf A}$
belong to $VB(2)$.
\end{proposition}

\begin{proof} The proof is based on Koszul complexes for the algebras
in question, which are defined as follows.

For an algebra $B$, let $B^e=B\otimes \bar B$.

Let $A={\mathcal A}$, ${\Bbb A}$, or ${\mathbf A}$.
Let $u=1\otimes x- x \otimes 1$ for ${\mathcal A}$,
and $u=X\otimes X^{-1}-1$ for ${\Bbb A}$ and ${\mathbf A}$.
Let $w=1\otimes p- p\otimes 1$ for ${\mathcal A}$ and ${\Bbb A}$, and
$w=P\otimes P^{-1}-1$ for ${\mathbf A}$.

We have an isomorphism of $A^e$-modules $A=A^e/I$, where $I=(u,w)$
is  an $A$-subbimodule generated by $u$ and $w$. As $u$ and $w$
commute, we can define the Koszul complex $K$ of $A^e$-modules:
$$
0\to {A}^e\otimes\Lambda^2 V\stackrel{d_1}{\to} {A}^e\otimes V
\stackrel{d_0}{\to}{A}^e\stackrel{\mu}{\to} {A},$$
where $V=\langle e_1,e_2\rangle\simeq \mathbb C^2$,
$d_1(a\otimes e_1\wedge e_2)=aw\otimes e_1-au\otimes e_2$,
$d_0(a\otimes e_1)=au$, $d_0(a\otimes e_2)=a w$ and $\mu$ is the
multiplication.

It is easy to check that the Koszul complex
is exact.

From the exactness of the Koszul complex
it follows that $A$ has Hochschild dimension 2.
Further, applying the functors
$Hom_{{A}}(\cdot,{A}^e)$ to the
truncated Koszul complex $K$, we get  the  complex:
\begin{gather*}
0\to \Hom_{{A}^e}({A}^e,{A}^e)\stackrel{d_0^*}{\to}
\Hom_{{A}^e}({A}^e,{A}^e)\otimes V^*
\stackrel{d_1^*}{\to}
\Hom_{{A}^e}({A}^e,{A}^e)\otimes \Lambda^2
V^*\to 0.
\end{gather*}
The homology of this complex is exactly $\Ext^*_{A^e}(A,A^e)$. But this
complex coincides with the Koszul complex $K$, because
$\Hom_{{A}^e}({A}^e,{A}^e)={A}^e$, and
the maps $d_0$, $d_1$ are dual to each other with respect to the
natural pairing. Thus $\Ext^*_{A^e}(A,A^e)$ is concentrated in degree 2,
and $\Ext^2_{A^e}(A,A^e)=A$, as desired.
\end{proof}

\subsection{Proof of Theorem \ref{classicalgps}}

For type $A$, $W=S_n$, so the formulas follow from Corollary \ref{vdb},
Proposition \ref{VB(2)}, and Propositions \ref{p1},\ref{p2},\ref{p3}

For type $B$, $W=S_n\ltimes \Bbb Z_2^n$, so
the formulas follow from Corollary \ref{vdb},
Proposition \ref{VB(2)}, Proposition
\ref{GA}, and Proposition \ref{AxZ2}.

\section{Cherednik algebras}

Let $A$ be an algebra over a field $K$ of characteristic zero.
Assume that $A\in VB(d)$, where $d\ge 2$ is an even number,
and $HH^0(A)=K, HH^1(A)=0$.
Theorem \ref{homo} implies that
for $n>1$, we have
$$
\dim HH^2(A(n))=\dim HH^2(A)+m,
$$
where $m=0$ if $d>2$ and $m=1$ if $d=2$.
Thus if $d>2$, all deformations of $A(n)$ come from deformations of $A$.
On the other hand, when $d=2$ then there is an additional deformation.
More precisely, we have the following result.

\begin{theorem}\label{cher}
If $d=2$ then the universal formal deformation of $A(n)$
depends on $\dim HH^2(A)+1$ parameters. More precisely, this
deformaion is a topologically free
algebra $A(n,c,k)$ over the ring $K[[c,k]]$
of functions on the formal neighborhood of zero
in $HH^2(A)\oplus K$. The algebra $A(n,c,0)$ is equal to $A_c(n)$,
where $A_c$ is the universal deformation of $A$.
\end{theorem}

\begin{proof} We have $HH^3(A)=HH_{-1}(A)=0$. Hence by Theorem \ref{homo}, \linebreak
$HH^3(A(n))=0$. Thus the deformations of $A(n)$ are unobstructed, and the
theorem follows.
\end{proof}

\begin{corollary}
Let $X$ be an affine symplectic surface
such that $H^1(X,\Bbb C)=0$.
Let $A_+$ be a quantization of $X$, $A=A_+[\hbar^{-1}]$.
Then the universal formal deformation of $A(n)$
depends on $\dim HH^2(X)+1$ parameters. More precisely, this
deformation is a topologically free
algebra $A(n,c,k)$ over the ring $K[[c,k]]$ ($K=\Bbb C((\hbar))$)
of functions on the formal neighborhood of zero
in $H^2(X,K)\oplus K$. The algebra $A(n,c,0)$ is equal to $A_c(n)$,
where $A_c$ is the universal deformation of $A$.
\end{corollary}

The algebra $A(n,c,k)$ is very interesting.
This is demonstated by the following examples.

1. Let $A={\mathcal A}$ be the Weyl algebra of rank 1.
Then $HH^2(A)=0$, and $A(n,c,k)=A(n,k)$ is the rational Cherednik algebra
of type $A_{n-1}$ (see e.g. \cite{EG}). This algebra is generated
by $x_1,...,x_n,p_1,...,p_n$ and $S_n$ with defining
relations
\begin{gather*}
\sigma x_i=x_{\sigma(i)}\sigma,\
\sigma p_i=p_{\sigma(i)}\sigma, \sigma\in S_n;\\
[x_i,x_j]=[p_i,p_j]=0,\\
[x_i,p_j]=ks_{ij}, i\ne j,\\
[x_i,p_i]=1-\sum_{j\ne i}ks_{ij},
\end{gather*}
where $s_{ij}$ is the permutation of $i$ and $j$.
It controls multivariable Bessel functions of type
$A$.

2. More generally, let $\Gamma$ be a finite subgroup of $SL_2(\Bbb C)$.
Let $A=\Bbb C[\Gamma]\ltimes {\mathcal A}$.
Then $\dim HH^2(A)=\nu(\Gamma)-1$, where $\nu(\Gamma)$ is the number
of conjugacy classes of $\Gamma$. Thus the universal deformation
$A(n,c,k)$ depends of $\nu(\Gamma)$ parameters. This is the
symplectic reflection algebra for the
wreath product $S_n\ltimes \Gamma^n$ defined in \cite{EG}.
For $\Gamma=\Bbb Z_2$ this algebra controls multivariable Bessel
functions of type $B$.

3. Let $A=\Bbb C[\Bbb Z_2]\ltimes B$,
where $B={\mathcal A},{\Bbb A}$, or ${\mathbf A}$.
Then $A(n,c,k)$ is the rational Cherednik algebra, degenerate
double affine Hecke algebra, and double
affine Hecke algebra of type $B_n$ \cite{Sa},
respectively. The symbol $c$ involves 1,2, and 4 parameters, respectively.
These algebra control Koornwinder polynomials and their
degenerations.

\begin{remarks} 1. Consider the more general case, where $HH^1(A)$ is
not necessarily zero. If $d=2$, Theorem \ref{homo} implies that
$HH^2(A(n))=HH^2(A)\oplus \Lambda^2 HH^1(A)\oplus K$.
We expect that the deformations of $A(n)$ in this case are still
unobstructed, and thus there exists a universal deformation
$A(n,c,\pi,k)$ depending on \linebreak $b_2+b_1(b_1-1)/2+1$ parameters (where
$b_i=\dim HH^i(A)$ and $c$, $\pi$, $k$ are the parameters of
$HH^2(A), \Lambda^2HH^1(A)$, and $K$, respectively),
such that $A(n,c,0,0)=A_c(n)$.

This expectation is the case for $A={\Bbb A},{\mathbf A}$. In this case,
$A(n,c,\pi,k)$ is the degenerate double affine Hecke algebra, respectively,
double affine Hecke algebra of type $A_{n-1}$ introduced by Cherednik.
These algebras control Jack and Macdonald polynomials of type
$A_{n-1}$, respectively.

2. For a general affine symplectic surface $X$,
one can only expect the existence of a {\it formal} deformation
$A(n,c,\pi,k)$. However, consider the special case when
$X$ admits a compactification $\bar X$ to which the Poisson bracket on $X$
extends by zero at the divisor at infinity.
In this case, motivated by the works \cite{Ar,EG,Ko2,NeSt}
and others, we expect the following somewhat vague picture
(for simplicity we restrict to the case $H^1(X,\Bbb C)=0$).

a) The deformation $A(n,c,k)$ is obtained by localization and completion
of an algebra $A(n,\hbar,c,k)$ over $\Bbb C$, where all
parameters are complex. In particular, the formal algebra $A_c=A(1,c)$
is obtained from the complex algebra $A(1,\hbar,c)$.

b) The algebra $A(n,0,c,k)$ is finite over its center.
For generic parameters, the spectrum $C(n,0,c,k)$
of the center $Z(n,0,c,k)$ of this algebra
is a smooth symplectic variety, which is a symplectic deformation of the
Hilbert scheme $Hilb_n(X)$ (it may be called the Calogero-Moser variety of
$X$). The algebra $A(n,0,c,k)$ is
the endomorphism algebra of a vector bundle of rank $n!$
on this variety.

c) A quantization of the variety $C(n,0,c,k)$ is provided
by the ``spherical subalgebra'' $eA(n,\hbar,c,k)e$,
where $e$ is the idempotent in $A(n,\hbar,c,k)$
which deforms the symmetrization idempotent of $\Bbb C[S_n]$ sitting in
$A(n,\hbar,c,0)$. In particular, the map $z\to ze$ defines
an isomorphism $Z(n,0,c,k)\to eA(n,0,c,k)e$
(``the Satake isomorphism'').

d) The variety $C(n,0,c,k)$ parametrizes
(in the sense of Berest-Wilson, \cite{BW}) stably free
ideals in $A(1,\hbar,c)$ (regarded as left modules).

3. We expect that in many interesting cases the algebras
$A(n,\hbar,c,k)$ can be described explicitly
(say, by generators and relations). It would be interesting,
for example, to find
such a description in the case when $X$ is $\Bbb CP^2$ with an elliptic
curve removed. In this case $A(1,\hbar,c)$ is obtained from the the
Odesskii-Artin-Tate elliptic algebra $E$ with three generators
see \cite{O}) by the formula
$A=E(1/K)_0$, where $K$ is the central element in $E$ of degree 3,
and subscript $0$ denotes elements of degree $0$.
\end{remarks}

\end{document}